\documentclass[12pt]{article}
\usepackage{latexsym,amsmath,amssymb,, amsthm}

\newif\ifdviwin

\usepackage[latin1]{inputenc}
\usepackage[english]{babel}
\usepackage{indentfirst}
\usepackage[mathscr]{eucal}
\usepackage{amssymb,amsmath,amsfonts}
\usepackage{fancybox,fancyhdr}
\usepackage{graphicx}

\newif\ifdviwin

\dviwintrue

\let\8=\infty \let\0=\emptyset

\def\flecha{\rightarrow}

\def\cte.{\mathop{\rm cte.}\nolimits}

\def\I{\mathop{\rm i }\nolimits}
\def\Re{\mathop{\rm Re }\nolimits}

\def\A{\mathbb{A}}

\def\B{\mathbb{B}}

\def\R{\mathbb{R}}

\def\C{\mathbb{C}}
\def\D{\mathbb{D}}

\def\H{\mathbb{H}}
\def\S{\mathbb{S}}

 \newtheorem{defi}{Definition}[section] \newtheorem{teo}[defi]{Theorem}
 \newtheorem{pro}[defi]{Proposition}
 \newtheorem{cor}[defi]{Corollary}

 \newtheorem{remark}{Remark}

\numberwithin{equation}{section}

\begin{document}
\mbox{}\vspace{0.4cm}\mbox{}

\begin{center}
\rule{15.2cm}{1.5pt}\vspace{0.5cm}

{\Large \bf Complete Flat Surfaces with two Isolated Singularities in Hyperbolic 3-space \footnote{Research partially supported by Ministerio de Educación y Ciencia
  Grant No. MTM2007-65249, Junta de Andalucía Grants No. FQM325, No. P06-FQM-01642 and PROCAD }}\\
\vspace{0.5cm} {\large Armando V. Corro $\mbox{}^a$, Antonio
Martínez$\mbox{}^b$ and Francisco Mil\'{a}n$\mbox{}^c$}\\
\vspace{0.3cm} \rule{15.2cm}{1.5pt}
\end{center}
  \vspace{1cm}
$\mbox{}^a$ Instituto de Matemática e Estatística,
Universidade Federal de Goiás, CP 131, Campus II, Samambaia, CEP 74.001-970-Goiâ, Brasil. \\
e-mail: corro@mat.ufg.br\vspace{0.2cm}

 \noindent $\mbox{}^b$, $\mbox{}^c$
Departamento de Geometría y Topología, Universidad de Granada,
E-18071 Granada, Spain. \\ e-mail: amartine@ugr.es ; milan@ugr.es
\vspace{0.2cm}

\noindent Date: February 27, 2009 
\vspace{0.3cm}

 \begin{abstract}
We construct examples of flat surfaces in $\H^3$ which are graphs over  a two-punctured horosphere and  classify complete embedded flat surfaces in $\H^3$ with only one end and  at most
two isolated singularities.
 \end{abstract}
\section{Introduction}
The theory of flat surfaces in $\H^3$  has undergone an important
development in the last few years. The starting point of this renewed
interest  has been the discovery in  \cite{GMM} that flat surfaces in
$\H^3$ admit a Weierstrass representation formula in terms of meromorphic
data, like the classical one for minimal surfaces in $\R^3$. This has
generated a great interest in such class of surfaces, even though the only
complete examples are the horospheres and hyperbolic cylinders (see
\cite{S}).

The last mentioned lack of complete examples has motivated an   important
advance in the problem of studying the singularities in these surfaces.
Questions such as their generic behaviour or the existence of complete
examples with singularities have been solved thanks to the works
 \cite{KUY}, \cite{KRSUY} and \cite{R}.

Contrarily to the minimal case,  flat surfaces in $\H^3$ can have isolated
singularities around which the surface is regularly embedded. Geometrically,  isolated singularities correspond to points where the Gauss map has not well defined limit.  Locally,
this kind of singularities have been classified  in \cite{GMi},  where is
proved that the class of flat surfaces that have $p\in \H^3$ as an
embedded isolated singularity  admits a one-to-one correspondence with the
class of analytic regular convex Jordan curves in the $2-$sphere.
But there are many interesting  problems in this theory that remain unsolved.  For example, we can quote the existence of compact or complete examples
with a finite number of isolated singularities. In this sense and up to
now,  the only known example of complete flat surface with isolated
singularities is the revolution one (also call \emph{the half hourglass})
which is a  graph over a  horosphere with only one point removed.
 The goal in this paper is to contribute to the understanding of this family of surfaces.

 The paper is organized as follows:

 Section 2 starts with some information about how flat surfaces in $\H^3$ can be represented by holomorphic data. Then, we deals with  the global behaviour of complete embedded flat surfaces with a finite number of isolated singularities, proving that any such surface is globally convex and, in particular, if it has only one end, then it is a graph over a finitely punctured horosphere.

 Section 3 is devoted to the construction of complete embedded surfaces with only two isolated singularities  and one end. The construction relies on the conformal representation of flat surfaces in $\H^3$ and the existence of conformal equivalences between a one punctured annulus and a horizontal slit domain in $\C$.

 Finally, in Section 4, we classify complete embedded flat surfaces in $\H^3$ with either one or two isolated singularities and only one end.

\section{Flat surfaces in $\H^3$ with isolated singularities}

 We consider the \emph{half-space} model of $\H^3$, that is,  $\H^3=\{(x_1,x_2,x_3)\in \R^3 : x_3>0\}$  endowed with metric
\begin{equation}
\langle,\rangle := \frac{1}{x_3^2} \left( dx_1^2 + dx_2^2 + dx_3^2\right),\label{im}
\end{equation}
of constant curvature $-1$
and with ideal boundary $\C_\infty = \{(x_1,x_2,0) : x_1,x_2\in \R\} \cup \{\infty\}$.

Let $\Sigma$ be a $2-$manifold and $\psi:\Sigma \longrightarrow \H^3$ be a
flat immersion. Then, from the Gauss equation, the second fundamental form
$d\sigma^2$ is definite and so $\Sigma$ is orientable and it inherits a
canonical Riemann surface structure such that the second fundamental form
$d\sigma^2$ is hermitian. This canonical Riemann surface structure
provides a conformal representation for the immersion $\psi$ that let to
recover any flat surface in $\H^3$ in terms of holomorphic data (see
\cite{GMM} and \cite{KUY} for the details).

For any $p\in \Sigma$, there exist $g(p),g_*(p)\in \C_\infty$ distinct
points in the ideal boundary such that the oriented normal geodesic at $\psi(p)$ is the geodesic in $\H^3$ starting
from $g_*(p)$ towards $g(p)$. The
maps $g,g_*:\Sigma \longrightarrow \C_\infty$ are called the
\emph{hyperbolic Gauss maps} and it is proved in \cite{GMM} that they are
holomorphic when we regard $\C_\infty$ as the Riemann sphere.

Kokubu, Umehara and Yamada investigated how to recover flat immersions with some admissible singularities (flat fronts) in terms of the hyperbolic Gauss maps. Making suitable the Theorem
2.11 in \cite{KUY} to the upper half-space model, we have
\begin{teo}[\cite{KUY}] \label{teokuy}
Let $g$ and $g_*$ be non-constant meromorphic functions on a Riemann surface $\Sigma$
such that $g(p)\neq g_*(p)$ for all $p\in \Sigma$. Assume that
\begin{enumerate}
\item all the poles of the 1-form $\frac{dg}{g-g_*}$ are of order 1, and
    \item $\Re\int_\gamma \frac{dg}{g-g_*}=0$, for each loop $\gamma$ on $\Sigma$.
\end{enumerate}
Set \begin{equation}
\xi:= c \exp \int \frac{dg}{g-g_*}, \qquad c\in\C\setminus \{0\}. \label{xi}
\end{equation}
Then,  the map $\psi=(\psi_1,\psi_2,\psi_3): \Sigma\longrightarrow \H^3$ given by
\begin{equation}
\psi_1 + \I \psi_2 = g - \frac{|\xi|^4 (g - g_*)}{|\xi|^4 + |g - g_*|^2}, \qquad \psi_3= \frac{|\xi|^2 |g - g_*|^2}{|\xi|^4 + |g - g_*|^2}\label{imhgm}
\end{equation}
is singly valued on $\Sigma$. Moreover, $\psi$ is a flat front if and only
if $g$ and $g_*$ have no common branch points.

Conversely any non-totally umbilical flat front can be constructed in this
way.
\end{teo}
Using (\ref{xi}) and (\ref{imhgm}), we see that for recovering $\psi$ we
only need a meromorphic function $g$ on $\Sigma$ and a harmonic function
$u:\Sigma\backslash \mathcal{P}_g\longrightarrow \R $, where
$\mathcal{P}_g$ is the set of poles of $g$,
\begin{equation}
u := \Re \int  \frac{dg}{g-g_*}.\label{harmonic}
\end{equation}
Moreover, the conditions in Theorem \ref{teokuy} say that $u$ and $g$
satisfy:
\begin{itemize}
\item[\textbf{(A)}] For each $p\in \mathcal{P}_g$, there exists a local coordinate $z$ vanishing at $p$ such  that $$ u + (b_g(p)+1)\log|z|$$ is harmonic in a neighborhood of $p$, where $b_g(p)$ is the branch number of $g$ at $p$.
\item[\textbf{(B)}] There exists a well-defined holomorphic function $F$ on $\Sigma$ such that
\begin{enumerate}
\item $du + \I \ast  du = F dg,$ where $\ast$ denotes the standard \emph{conjugation} operator acting on 1-forms, and
\item $g$ and $g - 1/F$ have no common branch points.
\end{enumerate}
\end{itemize}
A straight forward computation let us to prove that the induced metric and
the second fundamental form of $\psi$ are given, respectively, by
\begin{equation}
ds^2 = \exp(-4u)|\exp(4u) (dF + F^2dg) - \overline{dg} |^2, \label{pff}
\end{equation}
\begin{equation}
d\sigma^2 = \exp(-4u) |dg|^2 - \exp(4u) |dF + F^2dg|^2. \label{sff}
\end{equation}

All these facts let us to obtain the following conformal representation:
\begin{teo}\label{cr}
Let $g$  be a non-constant meromorphic function on a Riemann surface $\Sigma$ with set of poles $\mathcal{P}_g$ and $u:\Sigma\backslash \mathcal{P}_g \longrightarrow \R $  a harmonic function satisfying \textbf{(A)} and  \textbf{(B)}. If (\ref{pff}) (or equivalently, (\ref{sff})) is a riemannian metric, then,  the map $\psi=(\psi_1,\psi_2,\psi_3): \Sigma\longrightarrow \H^3$ given by
\begin{equation}
\psi_1 + \I \psi_2 = g - \psi_3 \exp(2u)\overline{F}, \qquad \psi_3=
\frac{ \exp(2u)}{1 + \exp(4u)|F|^2}, \label{immersion}
\end{equation}
is a well-defined  flat immersion.

Conversely, any flat immersion in $\H^3$ can be constructed in this manner.
\end{teo}
\begin{remark}
The hyperbolic Gauss map $g$ defines a horosphere congruences, for which
$\psi(\Sigma)$ and $g(\Sigma)$ are envelopes and $2 \exp(2u)$ is the
radius function of each horosphere (see \cite{C}).
\end{remark}
\begin{defi} Let $\Sigma$ be a differentiable surface without boundary,
 $\psi: \Sigma\flecha \H^3$  a continuous map  and $\mathcal{F}= \{p_1,\cdots,p_n\}\subset \Sigma$ a finite set. We say that $\psi$  is a complete flat immersion with isolated singularities $\psi(p_1), \cdots , \psi(p_n)$, if  $ \psi$ is a flat immersion in $\Sigma\setminus \mathcal{F}$ but  $\psi $ is not $C^1$ at the points $p_1,\cdots,p_n$, and every divergent curve in $\Sigma$ has infinite length for the induced (singular) metric.
\end{defi}
\begin{pro}\label{pro1}
Let  $\psi: \Sigma\flecha \H^3$ be a complete flat immersion with $\psi(\mathcal{F})$ as set of isolated singularities. Then there is a compact Riemannian surface $\overline{\Sigma}$, $n$ disjoint discs $\mathcal{D}_1, \cdots, \mathcal{D}_n \subset \overline{\Sigma}$ and finitely many points $q_1, \cdots, q_m \in \overline{\Sigma}\setminus \mathcal{D} $, where $\mathcal{D} = \mathcal{D}_1\cup \cdots \cup \mathcal{D}_n$ such that $\Sigma\setminus \mathcal{F}$ endowed with the conformal structure induced by the second fundamental form has the conformal type of $\overline{\Sigma}\setminus \{\{q_1,\cdots,q_m\}\cup \mathcal{D}\}$

The points $q_1,\cdots,q_m$ are called the ends of $\psi$.
\end{pro}
\begin{proof}
Let $\mathcal{K}\subset \Sigma$ be a closed disk containing $\mathcal{F}$ in its interior. From (\ref{pff}) and (\ref{sff}), we have that $ds^2\leq 2 \exp (-4 u) |dg|^2 $ in $\Sigma \setminus \overset{\circ}{\mathcal{K}}$. Thus the flat metric $\exp (-4 u) |dg|^2$ is complete and it follows from a classical result of Huber, \cite{H}, and Osserman, \cite{O}, that $\Sigma \setminus \overset{\circ}{\mathcal{K}}$ is conformally a compact Riemann surface with compact boundary and finitely many points $\{q_1,\cdots,q_m\}$ removed. Then the proposition follows because around an isolated singularity we have the conformal structure of an annulus (see Section 5 in \cite{GMi}).
\end{proof}
From Lemma 1 in \cite{GMM}, we know that a complete flat end in $\H^3$ must be conformally to a punctured disc.  Then, the following assertion follows as in \cite{Yu}, (see also the Appendix for details).
\begin{pro}\label{pro2}
Each embedded complete end of a flat surface in $\H^3$  is biholomorphic to a punctured disc and the hyperbolic Gauss map $g$ extends meromorphically to the punctured, that is, the end must be regular.
\end{pro}
\begin{teo}\label{teo1}
 If $\psi: \Sigma\flecha \H^3$ is a complete flat embedding with $\psi(\mathcal{F})$ as  set of isolated singularities, then $\psi$ is globally convex.
\end{teo}
\begin{proof}
   Consider  the Klein model for $\H^3$, that is, the diffeomorphism from $\H^3$ into the open
unit ball $\B^3 \subset \R^3$ given by $ \mathfrak{K}: \H^3 \flecha  \B^3$,
$$
\mathfrak{K}(y_1,y_2,y_3) =\left( \frac{2y_1}{\|y\|^2+1},\frac{2 y_2}{\|y\|^2+1},\frac{\|y\|^2-1}{\|y\|^2+1} \right),$$
 for any $y=(y_1,y_2,y_3)\in \H^3$ and where by $\| .\|$ we denote the usual Euclidean norm.

This map is totally geodesic, and thus, it  preserves convexity. In particular, flat surfaces in $\H^3$ are mapped into convex surfaces in $\R^3$. Moreover, the ideal boundary of $\H^3$ is mapped via $\mathfrak{K}$ to the the unit 2-sphere $\S^2$ of $\R^3$.

From the above Propositions, we have that $\mathfrak{K}(\psi(\Sigma))$ is a compact locally convex surface in $\R^3$ with a finite number of peaks which correspond to the ends and the isolated singularities. Moreover, from Theorem 11 in \cite{GMi} and having in mind that the ends are regular, we have that around each peak the surface is a convex graph over a plane passing through the peak.  Then, it is clear that  $\psi$ must be globally convex and we conclude the proof.
\end{proof}
\begin{cor}\label{coro}
Every  complete flat embedding $\psi:\Sigma \flecha \H^3$ with a finite number of  isolated singularities and only one  end is a graph over a finitely punctured horosphere.
\end{cor}

\section{Canonical examples}
In this Section we shall describe examples of complete flat embedding with  only one end and at most two isolated singularities.
\begin{figure}
\begin{center}
\DeclareGraphicsExtensions{eps}
\includegraphics[width=0.35\textwidth]{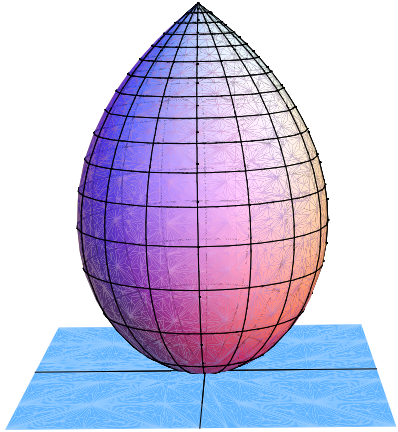}\\
\caption{Complete flat surface with only one isolated singularity}
\label{1singularidad}
\end{center}
\end{figure}
\subsection{Rotational examples}
It is known, see \cite{GMM} and \cite{KUY}, that a half hourglass is a revolution flat
complete embedding  with one isolated singularity
and one end,  which admits a conformal  parametrization  $\psi:\Sigma \longrightarrow \H^3$,   given by (\ref{imhgm}) with
$$g(z) = z, \;  \; g_*(z) = \frac{a + 1}{a - 1} z,$$
where $z \in \Sigma = \D^*_r = \{z \in \C \; / \; 0 < | z | < r \}$, $4 r^{2a} = 1 - a^2$ and  $a \in ]0, 1[$.

In this case, the   singularity is $\psi(\S_r) = (0, 0, b)$, with $\S_r =
\{z \in \C \; / \; | z | = r \}$ and $b \in \R$, the end is $\psi(0) = (0,
0, 0)$ and the function
$$ R(z) = F(z) g(z) = \frac{g(z)}{g(z) - g_*(z)}$$
is constant, see Figure \ref{1singularidad}.
\subsection{Examples with two isolated singularties}
We are going to construct examples of complete flat embedding with only two isolated singularities and one end.

 We know from Proposition \ref{pro1} and Corollary \ref{coro}, it admits a conformal parametrization $\psi:\A_r^\star\flecha \H^3$ of a punctured annulus $\A_r^\star=\A_r\setminus\{z_0\}$ in $\H^3$,  where
$$\A_r = \{z \in \C \; / \; r < | z | < 1 \},$$
 $0 < r < 1$ and $ z_0\in \A_r$ is the end.

Consider the annular Jacobi theta function
given by
\begin{eqnarray}
\vartheta_1(z) = C \left(1 - \frac{1}{z}\right) \prod_{k=1}^\infty (1 - r^{2k}z)(1 - r^{2k}/z), \label{theta}
\end{eqnarray}
here,  $C = \prod_{k=1}^\infty (1 - r^{2k})$. It is clear that it satisfies
\begin{eqnarray}
\vartheta_1(z) = \overline{\vartheta_1(\overline{z})} = - r^2 z \vartheta_1(r^2 z) = - \frac{1}{z} \vartheta_1(1/z), \; \; \vartheta_1(z/r^2) = - z \vartheta_1(z) \label{pr1}
\end{eqnarray}
and by deriving
\begin{eqnarray}
\vartheta_1'(z) & = & - r^2 \vartheta_1(r^2 z) - r^4 z \vartheta_1'(r^2 z) = \frac{1}{z^2} \vartheta_1(1/z) + \frac{1}{z^3} \vartheta_1'(1/z), \nonumber \\ \vartheta_1'(z/r^2) & = & - r^2 \vartheta_1(z)  - r^2 z \vartheta_1'(z). \label{pr2}
\end{eqnarray}
Thus, for any  $z_j \in ]-1, -r[$, one can see that the classical holomorphic bijection
$q_j : \A_r \backslash \{z_j\} \longrightarrow \C \backslash (I_1 \cup
I_r)$ given by
\begin{eqnarray}
q_j(z) = - \frac{\vartheta'_1(z_j/z)}{z \vartheta_1(z_j/z)} - \frac{z \vartheta'_1(z_j z)}{\vartheta_1(z_j z)}, \label{q}
\end{eqnarray}
maps the circles $\S_1$ and $ \S_r$, onto two real intervals $I_1$ and  $I_r$, respectively.

Actually, $q_j$ is characterized as the unique (up to real additive
constants) holomorphic map in $\A_r \backslash \{z_j\}$, which maps each
boundary component of $\A_r$ onto a real interval and has a simple pole of
residue $1$ at $z_j$, see \cite{A}.

Given  $z_0, z_1, z_2 \in ]-1, -r[$, we define the following  holomorphic
function $R : \A_r \backslash \{z_0\} \longrightarrow \C$
\begin{eqnarray}
R(z) = a q_0(z) + b, \label{R}
\end{eqnarray}
where $a$ and $b$ are real constants, determined by $R(z_1) = 1$, $R(z_2)
= 0$ and such that  $0<R<1$ on the boundary of $\A_r$, $\partial \A_r$. It is not a restriction to assume that
\begin{eqnarray}
R(\S_1) = [R(-1), R(1)] \subset \; ]0, 1[, \; R(\S_r) = [R(r), R(-r)] \subset \; ]0, 1[, \label{rs}
\end{eqnarray}
with $R(1) < R(r)$.

Then, $R'(\widetilde{z}) = 0$ only for $\widetilde{z} \in \{\pm 1, \pm r
\}$, and making an analysis of the set $q_0^{-1}(\R)$, one can see that
\begin{eqnarray}
R([-1, z_2] \cup \S_1 \cup [r, 1] \cup \S_r \cup [z_1, -r]) = [0, 1], \nonumber \\ \vspace{.3cm}
R(]z_2, z_0[) = ]-\infty, 0[, \; \; R(]z_0, z_1[) = ]1, +\infty[ \nonumber
\end{eqnarray}
and $z_0 \in \; ]z_2, z_1[$.

Moreover, by the above mentioned characterization of the holomorphic functions $q_j$
one has
\begin{eqnarray}
\frac{R'(z_1)}{R(z) - 1} = q_1(z) - c_1, \; \; \frac{R'(z_2)}{R(z)} = q_2(z) - c_2, \label{rp}
\end{eqnarray}
where $$c_1 = q_1(z_0) = q_1(z_2) + R'(z_1), \;  \; c_2 = q_2(z_0) =
q_2(z_1) - R'(z_2) \in \R.$$

With the above notations, we have:

\begin{pro}
If $z_0, z_1, z_2 \in ]-1, -r[$ and $m \in \R$ satisfy
\begin{itemize}
\item[\textbf{(C1)}] $m + c_1 z_1 - z_1 R'(z_1) = - z_2 R'(z_2),$
\item[\textbf{(C2)}] $c_1 z_1 - c_2 z_2 - 2 = 0,$
\item[\textbf{(C3)}] $z_1 z_2 r^{2(m+2)} = 1.$
\end{itemize}
Then the functions $g : \A_r \longrightarrow \C$ and $u : \A_r \backslash
\{z_0\} \longrightarrow \R$ given by
\begin{eqnarray}
g(z) = \sqrt{\frac{R(z)}{1 - R(z)} \frac{Q_1(z)}{Q_2(z)} z^{-2}}, \; \; \;
\label{g} u(z) = \frac{1}{2} \log{\left|\frac{Q_1(z)}{1 - R(z)} z^m
\right|} \label{u}
\end{eqnarray}
with
\begin{eqnarray}
Q_j(z) = \frac{\vartheta_1(z_j/z)}{\vartheta_1(z_jz)}, \; \; j = 1, 2, \label{Q}
\end{eqnarray}
satisfy the conditions of the Theorem \ref{cr}. Thus, $\psi: \A_r \setminus\{z_0\}\flecha \H^3 $  given by (\ref{immersion}) is a
well-defined flat surface with $\psi(\S_1)$ and $\psi(\S_r)$ as isolated singularities.
\end{pro}
\begin{proof}
As $z_j$ is a simple zero of $Q_j$, it is clear  that $g$ is a holomorphic
function, without zeros in $\A_r$, and $$u(z) - \frac{1}{2}\log|z-z_0|$$
is a harmonic function in $\A_r$. So, the condition \textbf{(A)} is
satisfied.

In order to check the condition \textbf{(B)}, we use that (\ref{q}) and
(\ref{Q}) give
$$\frac{d\log{Q_j(z)}}{dz} = \frac{z_j}{z} q_j(z).$$
Thus, from (\ref{rp}), \textbf{(C1)}, \textbf{(C2)} and (\ref{u}) we get
\begin{eqnarray}
2 F(z) g'(z) & = &  \frac{R'(z)}{1 - R(z)} + \frac{z_1}{z} q_1(z) +
\frac{m}{z} \nonumber \\ &=& \frac{R'(z)}{1 - R(z)} + \frac{z_1}{z}  \frac{R(z) R'(z_1)}{R(z) - 1} + \frac{m + c_1 z_1 - z_1 R'(z_1)}{z}
\end{eqnarray}
and
\begin{eqnarray}
\frac{2 g'(z)}{g(z)} & =&  \frac{R'(z)}{R(z)(1 - R(z))} + \frac{z_1}{z} q_1(z) - \frac{z_2}{z} q_2(z) - \frac{2}{z} \nonumber \\ & =&   \frac{R'(z)}{R(z)(1 - R(z))}
+ \frac{z_1}{z}  \frac{R'(z_1)}{R(z) -  1} - \frac{z_2}{z}  \frac{R'(z_2)}{R(z)} + \frac{c_1 z_1 - c_2 z_2 - 2}{z} . \label{gp}
\end{eqnarray}
Thus, we  conclude there exists a holomorphic function $F : \A_r \backslash \{z_0\} \longrightarrow \C$  given by
\begin{eqnarray}
F(z) = \frac{R(z)}{g(z)} \label{F}
\end{eqnarray}
which satisfies the condition \textbf{(B)}. It is clear that $g$ and $g- 1/F$ have no common branch points, because $R'\neq 0$ in $\A_r$.

On the other hand, if $z \in \S_1$, then from (\ref{immersion}), (\ref{pr1}),
(\ref{rs}), (\ref{u}), (\ref{Q}) and (\ref{F}), we obtain
\begin{eqnarray}
exp(2u(z)) = \frac{1}{1 - R(z)}, \; \; |g(z)|^2 = \frac{R(z)}{1 -
R(z)} \label{g1}
\end{eqnarray}
and so
\begin{eqnarray}
\psi(\S_1) = (0, 0, 1). \label{s1}
\end{eqnarray}
Similarly, if $z \in \S_r$, then
\begin{eqnarray}
exp(2u(z)) = \frac{|z_1| r^{m+1}}{1 - R(z)}, \; \; |g(z)|^2 =
\frac{R(z)}{1 - R(z)} \frac{z_1}{z_2} r^{-2} \label{gr}
\end{eqnarray}
and by using \textbf{(C3)}
\begin{eqnarray}
\psi(\S_r) = (0, 0, |z_1| r^{m+1}). \label{sr}
\end{eqnarray}
\end{proof}
We will call {\bf canonical examples} to those flat immersions obtained as in the above proposition.

\begin{figure}
\DeclareGraphicsExtensions{eps}
\begin{center}
\includegraphics[width=0.4\textwidth]{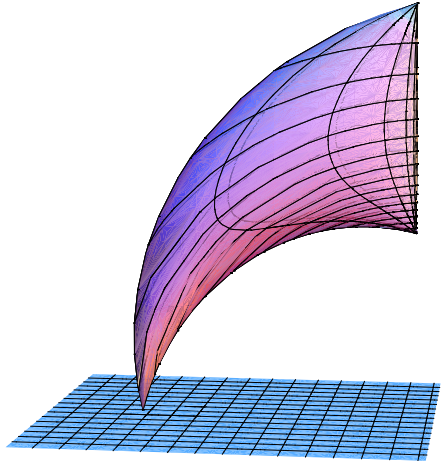}\\
\caption{Complete flat surface with two isolated singularities.}
\label{2picos}
\end{center}
\end{figure}

The following result proves that there exists a large family of canonical examples.
\begin{pro}
For any $r \in \; ]0, 1[$ and $s \in \; ]-1, 0[$, there exist $m \in \;  ]-3, -2[$ and $z_0, z_1, z_2 \in \; ]-1, -r[$, $z_2 < z_0 < z_1$, which satisfy the conditions \textbf{(C1), (C2), (C3)}, with $s = - z_2 c_2 = - z_2 q_2(z_0)$.

In particular, for $s = -1/2$, there is a solution with $m = -5/2$ and $z_0^2 = r = z_1 z_2$.
\end{pro}
\begin{proof}
From (\ref{q}) and (\ref{rp}), the condition \textbf{(C1)} can be written
\begin{eqnarray}
m = 2 z_1 z_2 \frac{\vartheta'_1(z_1 z_2)}{\vartheta_1(z_1 z_2)} - 1 + z_2 q_2(z_0) = 2h(z_1 z_2) - 1 - f_0(z_2),
\end{eqnarray}
where  $h$ and $f_0$ are the functions given by
\begin{eqnarray}
h(z) = \frac{z \vartheta'_1(z)}{\vartheta_1(z)}, \; \; f_0(z) = h(z/z_0) + h(zz_0). \label{hf}
\end{eqnarray}
 We are going to use some properties of $h$ and $f_0$. First, from (\ref{pr1}), (\ref{pr2}) and (\ref{hf}), the function $h$ verifies
\begin{eqnarray}
h(z) = 1 + h(r^2 z), \; \; h(z) + h(1/z) = - 1 \label{hz}
\end{eqnarray}
for any $z$. In particular $-1 = h(r) + h(1/r) = h(r) + 1 + h(r^2/r)$ gives
\begin{eqnarray}
h(r) =  - 1. \label{hr}
\end{eqnarray}
Moreover,  from (\ref{theta}), (\ref{hf}), (\ref{hz}) and (\ref{hr}), we obtain
\begin{eqnarray}
h(]r^2, r[) \; = \; ]-1, +\infty[, \; \; h(]r, 1[) \; = \; ]-\infty, -1[, \nonumber \\
f_0(]-1, z_0[) \; = \; ]-1, +\infty[, \; \; f_0(]z_0, -r[) \; = \; ]-\infty, -2[, \label{fi}
\end{eqnarray}
for any $z_0 \in \; ]-1, -r[$. See Figure \ref{hfunction} and Figure \ref{f0function}.
\begin{figure}
\DeclareGraphicsExtensions{eps}
\begin{center}
\includegraphics[width=0.3\textwidth]{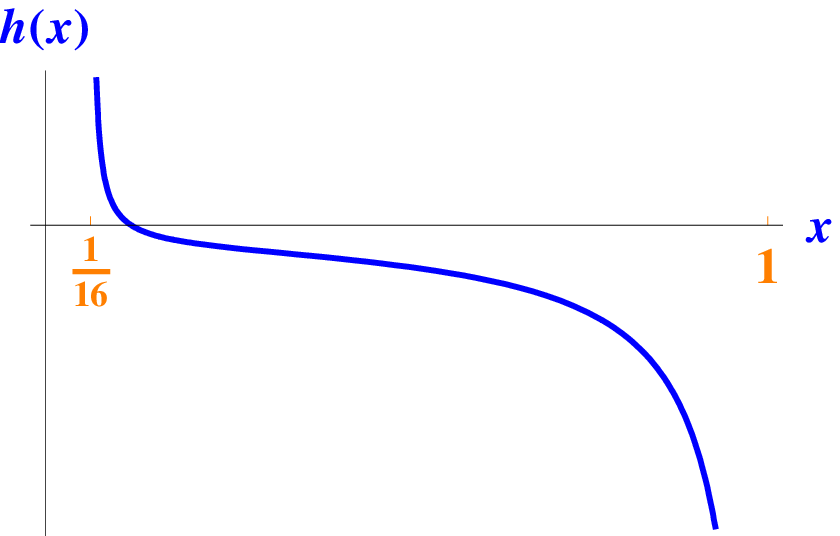}\\
\caption{Function $h(x)$, $x\in ]r^2,1[$ with $r=1/4$}
\label{hfunction}
\includegraphics[width=0.4\textwidth]{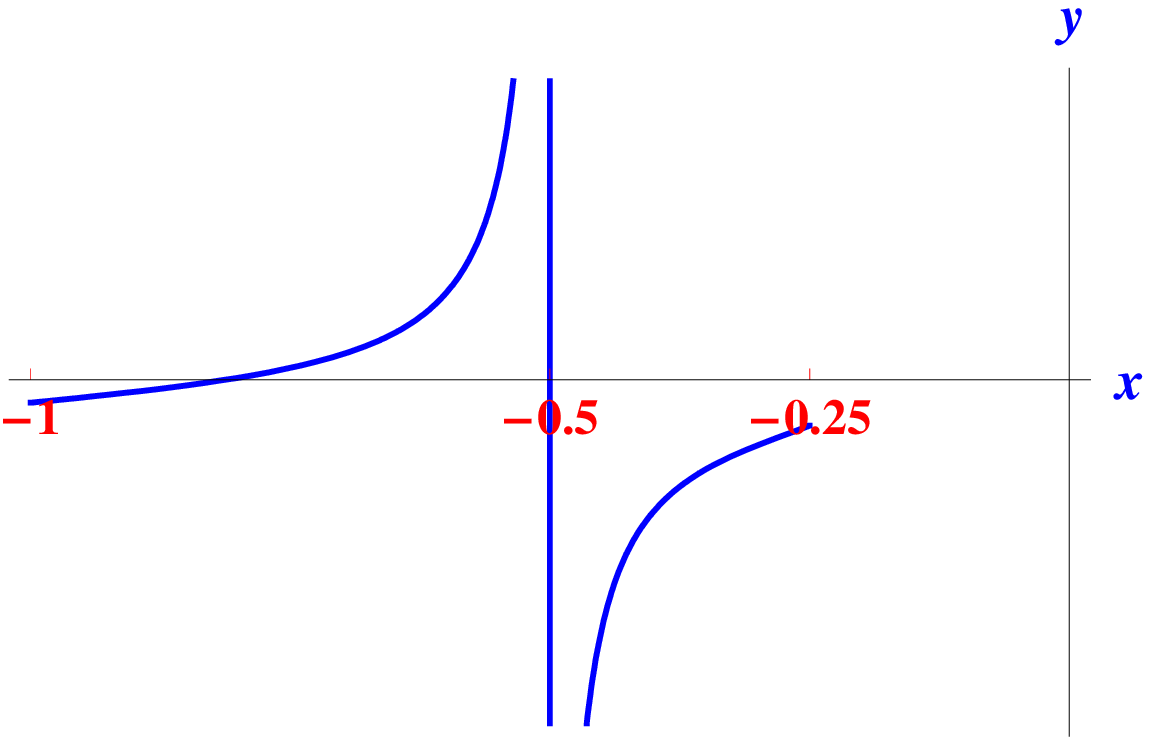}\\
\caption{Function $y=f_0(x)$, $x\in ]-1,-r[$ with $r=1/4$ and $z_0=-1/2$}
\label{f0function}
\end{center}
\end{figure}
But, the conditions \textbf{(C1), (C2), (C3)} are writing as
\begin{eqnarray}
m & = & 2h(r^{-2(m+2)}) - 1 - f_0(z_2), \nonumber \\ - 2 & = & f_0(z_1) - f_0(z_2), \label{cond} \\  z_1 z_2 & = & r^{-2(m+2)}, \nonumber
\end{eqnarray}
and from (\ref{fi}), if $s \in \; ]-1, 0[$, then there exits $m \in \;  ]-3, -2[$ such that
\begin{eqnarray}
m = 2 h(r^{-2(m+2)}) - 1 - s. \label{m}
\end{eqnarray}
In the same way,  for any $z_0 \in \; ]-1, -r^{-(m+2)} + \varepsilon[$, one gets $z_{20} \in \; ]-1, z_0[$ such that
\begin{eqnarray}
f_0(z_{20}) = s, \;  \;  z_{20} < z_0 <  \frac{r^{-2(m+2)}}{z_{20}} = z_{10},   \label{dos}
\end{eqnarray}
for an appropriate  $\varepsilon > 0$.

Finally, from (\ref{hf}), (\ref{cond}), (\ref{m}) and (\ref{dos}), we only need to find $z_0$ such that
\begin{eqnarray}
s- 2 = f_0(z_{10}) =  h(z_{10}/z_0) + h(z_{10} z_0) = \widetilde{f}(z_0).  \nonumber
\end{eqnarray}
This $z_0$ exists, because
\begin{eqnarray}
\lim_{z_0 \rightarrow -1} \widetilde{f}(z_0) = 2 h(r^{-2(m+2)}) = m + 1 + s > s - 2 \nonumber
\end{eqnarray}
and from (\ref{fi})
\begin{eqnarray}
\lim_{z_0 \rightarrow z_{10}} \widetilde{f}(z_0) = -\infty. \nonumber
\end{eqnarray}

In particular, for $s = - 1/2$, it is clear that (\ref{m}) is satisfied with $m = - 5/2$ and (\ref{hr}). Then, by taking $z_0 = - \sqrt{r}$, $z_2 \in \; ]-1, z_0[$ such that $f_0(z_2) = s$ and $z_1 z_2 = r = z_0^2$, (\ref{hf}) and (\ref{hz}) give
\begin{eqnarray}
f_0(z_1) = f_0\left(\frac{r}{z_2}\right) = h\left(\frac{r}{z_2 z_0}\right) + h\left(\frac{r z_0}{z_2}\right) = h\left(\frac{z_0}{z_2}\right) + h\left(\frac{r^2}{z_2 z_0}\right)  \nonumber
\\ = - 2 - h\left(\frac{z_2}{z_0}\right) + h\left(\frac{1}{z_2 z_0}\right) = - 3 - f_0(z_2) = -3 - s  = -2 + s \nonumber
\end{eqnarray}
and we finish the proof.
\end{proof}
\begin{teo}
Each  canonical example $\psi:\A_r \backslash \{z_0\} \longrightarrow \H^3$  is a complete  flat embedding with two isolated singularities and one end, (see Figure \ref{2picos}).
\end{teo}
\begin{proof}
From (\ref{pff}), (\ref{s1}) and (\ref{sr}), we have
\begin{eqnarray}
\exp(4u) \ (dF + F^2 dg ) = \overline{dg} \label{em}
\end{eqnarray}
on  $\partial \A_r$, and we deduce that $g'(z) \neq 0, \,  \forall z \in \partial \A_r$. Otherwise, if $g'(\widetilde{z}) = 0$, for some $\widetilde{z} \in \partial \A_r$, (\ref{F}) and (\ref{em}) imply $F'(\widetilde{z}) = 0 = R'(\widetilde{z})$, (that is, $\widetilde{z} \in \{\pm 1, \pm r \}$). But, in this case, (\ref{gp}) and \textbf{(C2)} give
\begin{eqnarray}
 0 = \widetilde{z} R(\widetilde{z}) (R(\widetilde{z}) - 1) \frac{2 g'(\widetilde{z})}{g(\widetilde{z})} =  z_1 R'(z_1) R(\widetilde{z}) + z_2 R'(z_2) (1 - R(\widetilde{z})) > 0, \nonumber
\end{eqnarray}
which is a contradiction
Then, from (\ref{rs}), (\ref{g1}) and (\ref{gr}), the holomorphic function $g$ is one to one on $\partial \A_r$, (because, it is a covering map with $g^{-1}(g(\widetilde{z})) = \{\widetilde{z} \}$, for $\widetilde{z} \in \{\pm 1, \pm r \}$), and using that  $\A_r \cup
\partial \A_r$ is compact, we deduce that the set
$$V = \{z \in \A_r  \; / \; g'(z) = 0\} = \{z \in \A_r \cup
\partial \A_r \; / \; g'(z) = 0\}$$ is finite. Again, as $g:\A_r \backslash g^{-1}(g(V)) \longrightarrow g(\A_r) \backslash g(V)$ is a covering map, (one to one on $\partial \A_r$), $g$ is a diffemorphism on $\A_r \backslash g^{-1}(g(V))$ and also on $\A_r$. Hence, $V = \emptyset$.

Now, from (\ref{u}) and (\ref{em}), the holomorphic function
$$p(z) = \left(\frac{Q_1(z)}{1 - R(z)} z^m \right)^2 \left(\frac{F'(z)}{g'(z)} + F^2(z) \right)$$
verifies $|p(z)| = 1$ on $\partial \A_r$. Then, by the maximum modulus principle, $|p(z)| < 1$ on $\A_r$ and the fundamental forms (\ref{pff}) and
(\ref{sff}) are positive definite and
$$ds^2 = \exp(-4u) |\overline{dg} - p \; dg|^2 \geq \exp(-4u) |dg|^2 (1 - |p|^2) = d\sigma^2.$$

As consequence, $\psi(\S_1)$ and $\psi(\S_r)$ are the unique singularities of the canonical examples. Moreover, outside a neighborhood of the singularities, we have $|p|^2 < 1 - \varepsilon$, for some $\varepsilon > 0$ and
$$ds^2 \geq \varepsilon \exp(-4u) |dg|^2 = \varepsilon \left| \frac{1 - R(z)}{Q_1(z) z^m} \right|^2 |dg|^2$$ is complete, because $R$ has a pole in $z_0$. As  $g$ is a diffemorphism, it follows, by (\ref{u}), that $\psi(z_0) = (g(z_0), 0)$ is the only end and it is embedded,  see \cite{GMM} and \cite{KUY}.

Finally, $\psi:\A_r  \longrightarrow \R^3$ induces a local diffemorphism, well
defined and continuous on the topological sphere obtained when $\S_1$
and $\S_r$ are identified with two points. So, $\psi$ is a covering map, with
embedded end, and we conclude that it is one to one.
\end{proof}

\section{Characterizations}
\begin{teo}
The revolution examples are the unique complete  flat embeddings in $\H^3$ with only one isolated singularity and one end.
\end{teo}
\begin{proof}
If
$\psi: \Sigma \longrightarrow \H^3$ is a complete  flat embedding with only one isolated singularity and one end,  we know from Proposition \ref{pro1} and Theorem \ref{teo1} that $\psi$ admits a conformal parametrization  $\psi: \D_r \backslash \{z_0\} \longrightarrow \H^3$
where $z_0 \in \; \D_r$ is the end and $\psi(\S_r)$ is the isolated singularity.

Now, up to isometries of $\H^3$, we have $\exp(2u(z_0)) = 0,$
\begin{eqnarray}
\psi(\S_r) = (0, 0, 1), \; \; \psi(z_0) = (g(z_0), 0) = (0, 0, 0). \label{rsf}
\end{eqnarray}
Hence, if $z \in \S_r$, (\ref{immersion}) and (\ref{rsf}) give
$$
g(z) = exp(2u(z)) \; \overline{F(z)}, \; \; exp(2u(z)) = 1 + |g(z)|^2
$$
and
$$
\overline{g(z) F(z)} = \frac{|g(z)|^2}{1 + |g(z)|^2}.
$$
That is, the holomorphic function $R: \D_r \backslash \{z_0\}
\longrightarrow \C$, defined by $R(z) = g(z) F(z)$, is real in $\S_r$ and
\begin{eqnarray}
|g(z)|^2 = \frac{R(z)}{1 - R(z)}, \; \; exp(2u(z)) = \frac{1}{1 - R(z)}. \label{rar1}
\end{eqnarray}
From (\ref{rsf}) we also obtain that the harmonic function $u: \D_r
\backslash \{z_0\} \longrightarrow \R$ is given by
$$u(z) = \ln|z - z_0|^n + \widetilde{u}(z)$$ with $n > 0$ and $\widetilde{u}:
\D_r \longrightarrow \R$ a harmonic function. Then
\begin{eqnarray}
du + \I \ast du = F dg = R \; \frac{dg}{g} \nonumber
\end{eqnarray}
has a simple polo in $z_0$ and, as $g(z_0) = 0$, (or $g(z_0) = \infty$),
$R$ is a holomorphic function on $\D_r$ and real on $\partial \D_r$, that
is, $R$ is a constant function, $R(z) = b$, $\forall z\in \D_r$ and
\begin{eqnarray}
exp(2u(z)) = a |g(z)|^{2b} \label{reu}
\end{eqnarray}
for any  $z \in \D_r$, where $b \in ]0, 1[$ and $a > 0$.

Consequently, from (\ref{immersion}), (\ref{rar1}) and (\ref{reu}), we conclude that $\psi$ is  the revolution example given by
\begin{eqnarray}
\psi(z) = \widetilde{\psi}(g) = \left(g \frac{1 - a^2 (b - b^2) |g|^{4b-2}}{1 + a^2 b^2 |g|^{4b-2}}, \frac{a |g|^{2b}}{1 + a^2 b^2 |g|^{4b-2}}  \right) \nonumber
\end{eqnarray}
with $g \in \D_s \backslash \{0\}$, for $s = \sqrt{\frac{b}{1-b}}$ and $a
= (1-b)^{b-1} b^{-b}$.

\end{proof}

\begin{teo}
Each  complete  flat embedding in $\H^3$  with only two isolated singularities and one end must be congruent to one of the canonical examples.
\end{teo}
\begin{proof}
If
$\psi: \Sigma \longrightarrow \H^3$ is a complete  flat embedding with only two isolated singularities and one end,  we have from Proposition \ref{pro1} and Theorem \ref{teo1} that $\psi$ admits a conformal parametrization  $\psi: \A_r \backslash \{z_0\} \longrightarrow \H^3$,
where $z_0 \in \; ]-1, -r[$ is the end and the singularities are the
points $\psi(\S_1)$ and $\psi(\S_r)$.

Also, up to isometries of $\H^3$, we can consider $\exp(2u(z_0)) = 0,$
\begin{eqnarray}
\psi(\S_1) = (0, 0, 1), \; \; \psi(\S_r) = (0, 0, c),  \label{sf}
\end{eqnarray}
with $c \in \R^+ \backslash \{1\}$.

Now, if $z \in \S_1$, (\ref{immersion}) and (\ref{sf}) give
\begin{eqnarray}
|g(z)|^2 = \frac{R(z)}{1 - R(z)}, \; \; exp(2u(z)) = \frac{1}{1 - R(z)} \label{ar1}
\end{eqnarray}
and if $z \in \S_r$, then
\begin{eqnarray}
|g(z)|^2 = \frac{c^2 R(z)}{1 - R(z)}, \; \; exp(2u(z)) = \frac{c}{1 - R(z)}, \label{ar2}
\end{eqnarray}
where $R: \A_r \backslash \{z_0\} \longrightarrow \C$ is the holomorphic
function $g F$.

Again $R$ is real on $\partial \A_r$ and
\begin{eqnarray}
du + \I \ast du = F dg = R \; \frac{dg}{g} \label{du}
\end{eqnarray}
has a simple pole in $z_0$.

However, in this case, $\log|g|^2$ is a harmonic function on $\A_r$,
because if $g$ has a zero or a pole in $z_0$, then $R(z) = b \in ]0, 1[$
and one obtains a revolution example with only one singularity. Moreover,
as $g'(z_0) \neq 0,$ since the end is embedded, $R$ has a simple polo in
$z_0$ and, by the characterization of $q_0$, must be (\ref{R}).

So, by the uniqueness of the Dirichlet problem for harmonic functions on
$\A_r$, (\ref{ar1}) and (\ref{ar2}) we get
\begin{eqnarray}
\log|g(z)|^2 = \log\left|\frac{R(z)}{1 - R(z)} \frac{Q_1(z)}{Q_2(z)}
z^{n}\right|, \nonumber
\end{eqnarray}
with $n \in \R$ such that
\begin{eqnarray}
c^2 = \frac{z_1}{z_2} r^n. \label{c2}
\end{eqnarray}
Thus, after a rotation, $g$ is given by (\ref{g}), (with $n = -2$, since
the end is embedded).

Finally, from (\ref{gp}) and (\ref{du}), we also have the harmonic
function $\widetilde{u}: \A_r \longrightarrow \R$,
\begin{eqnarray}
\widetilde{u}(z) = u(z) - \frac{a g'(z_0)}{g(z_0)} \log |z-z_0| = u(z) -
\frac{1}{2} \log |z-z_0|, \nonumber
\end{eqnarray}
determined by (\ref{ar1}) and (\ref{ar2}). These conditions coincide with
(\ref{g1}) and  (\ref{gr}) on $\partial \A_r$, if we take $m \in \R$ such
that
\begin{eqnarray}
c = |z_1| r^{m+1} \label{c3}
\end{eqnarray}
and so, the harmonic fuction $u: \A_r \backslash \{z_0\} \longrightarrow
\R$ is (\ref{u}).

Then, from (\ref{u}), (\ref{du}), (\ref{c2}) and (\ref{c3}) one deduces
\textbf{(C1)}, \textbf{(C2)}, \textbf{(C3)} and the flat surface is one of
the canonical examples.

\end{proof}
\begin{remark}
From the above proofs, it is clear that there are not compact embedded
flat surfaces, with less than three isolated singularities, because $R$ is
constant only for the revolution examples.
\end{remark}
\section{Appendix}
As we have remarked before Proposition 2.5, from \cite{GMM}, we know that a complete flat end in $\H^3$ must be conformally to a punctured disc $\D^*$ and admits a conformal parametrization $$\psi: \D^* \longrightarrow \H^3,$$ with associated Weierstrass data $(f(z), h(z)dz)$, such that
\begin{eqnarray}
 \frac{h'}{h} = \sum_{i=-1}^\infty p_i z^i, \; \; f h^2 =  \sum_{i=-m}^\infty q_i z^i,          \label{wd}
\end{eqnarray}
where $m \geq 3$ if and only if the end is irregular. Moreover, the hyperbolic Gauss maps are given by
\begin{eqnarray}
g(z) = \frac{X_2(z)}{X_1(z)}, \; \;  g_*(z) = \frac{X'_2(z)}{X'_1(z)}      \label{hg}
\end{eqnarray}
being  $X_1, X_2$ linearly independent solutions of the ordinary linear differential equation
\begin{eqnarray}
X'' -   \frac{h'}{h} X' - f h^2 X = 0.    \label{ode}
\end{eqnarray}
This equation is of the same type that the equation (E1) studied by Yu in \cite{Yu}. In particular, he proved that its fundamental solutions in a sector domain take the following forms
\begin{eqnarray}
X_1(z) = z^{a+\sigma} (1 + A(z)) \exp(- \zeta), \; \; X_2(z) = z^{a-\sigma} (1 + B(z)) \exp(\zeta)    \label{fs}
\end{eqnarray}
where $$\zeta = \frac{\sigma_{n-2}}{z} + ... + \frac{\sigma_0}{z^{n-1}},$$ $A, B$ are analytic functions, which tend to zero as $z$ tends to zero, and the numbers $a, \sigma, \sigma_{n-2}, ..., \sigma_0$ and $n$ depending of the coefficients (\ref{wd}) and $m \geq 3$.

Finally, from Lemma 3 in \cite{Yu}, one can choose a sector domain which contains an essential direction of the function $$z^{-\sigma} \exp(\zeta).$$ Then, by using (\ref{hg}), (\ref{fs}) and the proof of Theorem 6 in \cite{Yu}, one gets
\begin{teo}
No irregular ends of flat surfaces in $\H^3$ are embedded.
\end{teo}
\vspace{.5cm}
\noindent\textbf{Acknowledgements}\\
This work was started during the first author's visit to the Department of Geometry and Topology at the University of Granada in 2008. He would like to thank the members of the department for their hospitality.
\def\refname{References}

\end{document}